\documentclass[12pt,reqno]{amsart}

\usepackage{amssymb,delarray}
\usepackage{graphicx}

\usepackage{epsf,amsfonts,amscd,latexsym}

\usepackage[russian]{babel}

\def\leq{\leqslant}
\def\geq{\geqslant}

\def\edvo{\rule {6pt}{6pt}} 

\newtheorem{thm}{Theorem}
\newtheorem{lem}{Lemma}
\newtheorem{pro}{Assertion}

\newtheorem{cor}{Corollary}

\newtheorem{rem}{Remark}
{\catcode`\@=11
\gdef\n@te#1#2{\leavevmode\vadjust{%
 {\setbox\z@\hbox to\z@{\strut#1}%
  \setbox\z@\hbox{\raise\dp\strutbox\box\z@}\ht\z@=\z@\dp\z@=\z@%
  #2\box\z@}}}
\gdef\leftnote#1{\n@te{\hss#1\quad}{}}
\gdef\rightnote#1{\n@te{\quad\kern-\leftskip#1\hss}{\moveright\hsize}}
\gdef\?{\FN@\qumark}
\gdef\qumark{\ifx\next"\DN@"##1"{\leftnote{\rm##1}}\else
 \DN@{\leftnote{\rm??}}\fi{\rm??}\next@}}

\begin{document}
\baselineskip=14pt plus 2pt

\title[]{On the conjugacy problem in group $\bf F/{N_1\cap N_2}$.}
\author[O.V.~Kulikova]{O.V.~Kulikova}
\address{
\newline Bauman MSTU
\newline olga.kulikova@mail.ru }

\dedicatory{}\subjclass{}

\begin{abstract}
Let $N_1$ (resp., $N_2$) be the normal closure of a finite symmetrized set $R_1$ (resp., $R_2$) of a finitely generated free group $F = F(A)$.
It is well-known that if $R_i$ satisfies the condition $C(6)$, then the conjugacy problem is solvable in
$F/N_i$.  In the present paper we prove that if $R_1\cup R_2$ satisfies the condition $C(6)$ and the presentation
$\langle A\, \mid
R_1,R_2\rangle$ is atorical, then the conjugacy problem is solvable in $F/{N_1\cap N_2}$. In particular, if $R_1\cup R_2$ satisfies the condition $C(7)$ then the conjugacy problem is solvable in $F/{N_1\cap N_2}$.

Bibliography: 13 items.
\end{abstract}

\maketitle \setcounter{tocdepth}{2}

\def\st{{\sf st}}

\setcounter{section}{0}
 \section*{ Introduction.}
Let $F=F(A)$ be a free group generated by a finite alphabet $A$.
Let $N_1$ (resp., $N_2$) be the normal closure of non-empty finite
set $R_1$ (resp., $R_2$) of elements
of $F$. Assume that $R_i$ ($i = 1, 2$) is symmetrized, i.e., all elements of $R_i$ are cyclically reduced and
for any $r$ of $R_i$ all cyclic permutations of
$r$ and $r^{-1}$ also belong to $R_i$.

We will use the following notations. Denote graphic (letter-by-letter) equality of words $u,v\in F$ by $u\equiv v$.
Denote free equality by $u=v$. If words $u,v\in F$ present equal elements in a group $H$, we will write: $u=v$ in $H$.

If two words $u, v\in F$ are equal both in the group $F/{N_1}$ and
in the group $F/{N_2}$, then they are evidently equal in the group
$F/{N_1\cap N_2}$. It is natural to ask whether the conjugation of
words $u$ and $v$ in $F/{N_1\cap N_2}$ follows from their
conjugation both in $F/{N_1}$ and $F/{N_2}$? The answer is
obviously negative.  As an example showing that one can consider
the free group $F = F(a,b,c)$, the sets $R_1 = \{a^{\pm 1}\}$,
$R_2 = \{b^{\pm 1}\}$ and the words $u \equiv c^2ba$, $v \equiv cbca$.

The aim of this paper is to find out conditions on $R_1$ and $R_2$
such that the solvability of the conjugacy problem in $F/{N_1\cap
N_2}$ follows from the solvability of the conjugacy problem in
$F/{N_1}$ and $F/{N_2}$.

Note that this problem is naturally associated with subdirect products.
Indeed, one can consider $F/{N_1\cap N_2}$ as a subgroup of the direct product of
$F/{N_1}$ and $F/{N_2}$, and $F/{N_1\cap N_2}$ is a subdirect product of $F/{N_1}$ and $F/{N_2}$. Conversely, given a subdirect product $H$ of groups $G_1$ and $G_2$, there exist normal subgroups $N_1$ and $N_2$ of some free group $F$ such that $F/N_i\cong G_i$ (i=1,2) and $F/{N_1\cap N_2}\cong H$.

In turn,  subdirect products of two groups are closely associated to the fibre product construction in the category of groups. Recall that, associated to each pair of short exact sequences of groups  $1\rightarrow L_i\rightarrow G_i\stackrel{{\Psi}_i}{\rightarrow} Q\rightarrow 1,$
$i=1,2$, one has the {\it fibre product} $H = \{(x,y)\in G_1\times G_2 | {\Psi}_1(x) = {\Psi}_2(y)\}$. It is shown in \cite{bridson_miller} that a subgroup $H\leq G_1\times G_2$ is a subdirect product of $G_1$ and $G_2$ if and only if there is a group $Q$ and surjections ${\Psi}_i: G_i \rightarrow Q$ such that  $H$ is the fibre product of ${\Psi}_1$ and ${\Psi}_2$.

The question about the solvability of the conjugacy problem for subdirect products has been already considered for some groups (see, for example, \cite{bridson_miller, miller, boumslag}). Thus in the paper of C.F. Miller \cite{miller} there is an example of a fibre product in which $G_1 = G_2$ are non-abelian finitely generated free groups, $L_1 = L_2$, ${\Psi}_1 = {\Psi}_2$, $Q$ is a finitely presented group with undecidable word problem, and the conjugacy problem in $H$ is unsolvable. So the natural question, whether  the solvability of the conjugacy problem in $F/{N_1\cap N_2}$ always follows from the solvability of the conjugacy problem both in $F/{N_1}$ and in $F/{N_2}$, has negative answer.  Since $Q$ is isomorphic to $F/N_1N_2$, it follows from the example of C.F. Miller that  the solvability of the word problem in $F/N_1N_2$ is necessary for the solvability of the conjugacy problem in $F/{N_1\cap N_2}$.

To formulate the main result of the
paper, we recall the definitions of some geometrical
objects, called pictures. Pictures were introduced in \cite{igusa, rourke}. These objects are a very useful tool in combinatorial group theory, and can be used in a variety of different ways (see, for example,  \cite  {pride, gener_b_pr} and references in these papers).

Let $N$ be the normal closure of a symmetrized set $R$ of the free group $F(A)$.

{\it A picture} $P$ over the presentation $ \hat{G}=\langle A\,
\mid R \rangle$ on an oriented surface $T$ is a finite collection
of \,"vertices"\, $V_1, . . . , V_n \in T$, together with a finite
collection of simple pairwise disjoint connected oriented
"edges"\, $E_1, . . . , E_m \in T\setminus\{ \{V_1, ..., V_n\}\cup
\partial T\}$ labelled
by words of $F(A)$. But these edges need not all connect two
vertices. An edge may connect a vertex to a vertex (possibly
coincident), a vertex to $\partial T$, or $\partial T$ to
$\partial T$. Moreover, some edges need have no endpoints at all,
but be circles disjoint from the rest of $P$, such edges are
called edges-circles.

In the paper we will only ever consider such paths on $T$, each of which does't pass through any vertex  and intersects the edges of $P$ only finitely many times (moreover, if a path intersects an edge then it crosses it, and doesn't just touch it).
If
we travel along an oriented path $\gamma$ in the positive direction, we encounter a
succession of edges $E_{i_1},...,E_{i_k}$ labeled by
$g_{i_1},...,g_{i_k}$ respectively. These labels form the word
$g_{i_1}^{\varepsilon_{i_1}}\cdot...\cdot
g_{i_k}^{\varepsilon_{i_k}}$, where
${\varepsilon_{i_j}}\in\{1,-1\}$ is a local intersection index of
$E_{i_j}$ and $\gamma$. This word will be called {\it the word
along the path $\gamma$} (or {\it the label of $\gamma$}) and
denoted by $Lab^{+}( \gamma )$. The subword $g_{i_j}^{\varepsilon_{i_j}}$ ($j=1,...,k$) will be called {\it the contribution of $E_{i_j}$ in the label of
$\gamma$}. Travelling along $\gamma$ in the
negative direction gives the word ${{Lab}^{-}( \gamma )} \equiv
{Lab}^{+}( \gamma )^{-1}$.

If a path $\gamma$ is closed, consider a point $p$ on $\gamma$ not
belonging to any edge of $P$. The word along $\gamma$ read from
$p$ will be denoted by ${Lab_p}^{+}( \gamma )$ or by ${Lab_p}^{-}(
\gamma )$ (depending on the direction of travelling along
$\gamma$). Changing the disposition of $p$ we obtain the same word
up to cyclic permutation. We will denote the word along the path
$\gamma$ by $Lab( \gamma )$ when the disposition of $p$ and the
direction of reading will not be essential.

For each vertex $V$ of $P$ consider a circle $\Sigma$ of a small radius
with center at $V$  and a point $p\in \Sigma$ not lying on any edge of
$P$. The word ${Lab_p}^{+}( \Sigma )$ is called {\it the label of the
vertex} $V$. To complete the definition of the picture over the
presentation $ \hat{G}=\langle A\, \mid R \rangle$ on the surface
$T$ it remains to require that the labels of all vertices in $P$
belong to $R$.

Below we will consider pictures on a surface $T$, where $T$ is a
torus (torical pictures), an annulus (annulus pictures) or a disk
(planar pictures).

For a planar picture {\it the boundary label} of the picture is
the word ${Lab_{\bar{p}}}^{+}( \bar{\Sigma} )$, where $\bar{\Sigma}$ is a
circle near the boundary of the disk $T$ and $\bar{p}\in \bar{\Sigma}$
is a point not belonging to any edge.

The following result is well-known (use Theorem 11.1 \cite{olsh}
and dualise):
\begin{lem}
Let $W$ be a non-empty word on the alphabet $A$. Then $W$
represents the identity of the group $\hat{G}=F/N$ if and only if
there is a planar picture over the presentation $\langle A\, \mid
R \rangle$ of $\hat{G}$ with the boundary label $W$.
\end{lem}

A {\it dipole} is two distinct vertices $V_1$ and $V_2$ of $P$
connected by an edge $E$ if there exists a simple path $\psi$
joining points $p_1$ and $p_2$ on the circles $\Sigma_1$ and $\Sigma_2$
around these vertices, passing along $E$ and not crossing any edge
or vertex such that ${Lab_{p_1}}^{+}(\Sigma_1)={Lab_{p_2}}^{-}(\Sigma_2)$ in
$F$.

A presentation $ \hat{G}=\langle A\, \mid R \rangle$ is called
{\it atorical} (see, for example, \cite{olsh}) if every connected torical picture over $\hat{
G}=\langle A\, \mid R \rangle$ having at least one vertex contains
a dipole.

The following theorem $1$ will be proved in Section \ref{sec1}.

\begin{thm}Let $F$ be a free group generated by a finite alphabet $A$, $N_1$ (resp., $N_2$) be the normal closure of non-empty finite
symmetrized set $R_1$ (resp., $R_2$) of elements
of $F$.

Let the following conditions hold for the group $G_i=F/N_i$
($i=1,2$):

1.1.  The conjugacy problem is solvable in $G_i$.

1.2. In
$G_i$, there exists an algorithm allowing for a
reduced word $x\in F$, $x\neq 1$, to determine all
$z\in F$ such that $x\in \langle z \rangle$ in $G_i$, and the number of such distinct elements $z$ of $G_i$ is finite.

Let the following conditions hold for the group $G=F/N_1N_2$:

2.1.  The membership problem for a cyclic subgroup is solvable in
$G$.

2.2. The presentation $G=\langle A\, \mid R_1\cup R_2 \rangle$ is
atorical.

Then the conjugacy problem is solvable in $F/{N_1\cap N_2}$.
\end{thm}

Note that Condition 2.2 of Theorem 1 provides the equality $N_1\cap N_2 = [N_1, N_2]$ for disjoint $R_1$ and $R_2$ (see, for example, \cite{ratcl, ok1}).

Recall the definition of small cancelation conditions $C(k)$ (\cite{lind}), used in Theorem 2 below. A nontrivial freely reduced word $b$ in $F$ is called a {\it piece} with respect to $R$ if there exist two distinct elements $r_1$ and $r_2$ in $R$ that both have $b$ as maximal initial segment, i.e. $r_1\equiv bc_1$ and $r_2\equiv bc_2$. Let $k$ be a positive integer.
$R$ is said to satisfy {\it the small cancelation condition $C(k)$}, if no element of $R$ can be written as a reduced product of fewer than $k$ pieces.

Using the notations of Theorem 1, we have:

\begin{thm}
If $R_1\cup R_2$ is a set satisfying the condition $C(6)$ and the
presentation $G=\langle A\, \mid R_1\cup R_2 \rangle$ is atorical,
then  the conjugacy problem is solvable in $F/{N_1\cap N_2}$.
\end{thm}

\noindent{\it Proof of Theorem $2$.} Let us show that Theorem $2$
follows from Theorem $1$. Since $R_1\cup R_2$ satisfies the
condition $C(6)$, the subsets $R_1$ and $R_2$ also satisfy the
condition $C(6)$. Therefore
Condition 1.1 for $G_1$ and $G_2$ follows from Theorem 7.6
\cite{lind}; Condition 2.1 follows from Theorem 1 \cite{bezv}; Condition 1.2 can be deduced from Theorem 1 \cite{bezv}, Theorem 2
\cite{bezv} and (if there is an element of finite order in $G_i$) Theorem 1.4 \cite{b_pr} with regard to Theorem 13.3 \cite{olsh}. \edvo

It is well known that the condition $C(7)$ is sufficient for
atoricity (the proof of it is similar to Theorem 13.3
\cite{olsh}). So by Theorem $2$ (using the notations of Theorem 1) we have the following:

\begin{cor}
If $R_1\cup R_2$ satisfies the condition $C(7)$, then  the
conjugacy problem is solvable in $F/{N_1\cap N_2}$.
\end{cor}

The author is grateful to N.V.Bezverkhnii and A. Muranov for useful conversations
during preparation of this paper and, most particular, A. Minasyan, who showed the author the example of C.F. Miller and the relationship between the present work and
subdirect and fibre products.

\section{Deduction of Theorem $1$ from  Assertion
$1$.}\label{sec1}

Below, no mentioning it explicitly, we will use the fact that Condition 1.1 (resp., 2.1) of Theorem 1 leads to the solvability of the word problem in $G_i$, $i\in\{1,2\}$ (resp., $G$).

Let $u$ and $v$ be two reduced words of $F$. For each
$i\in\{1,2\}$ by Condition $1.1$ of Theorem $1$ there exists an algorithm which decides whether $u$ and $v$ present
conjugated elements in $G_i=F/N_i$. If $u$ and $v$ turn out to be
not conjugated in $G_i$ for at least one of the $i$, then $u$ and $v$
are not conjugated in $F/{N_1\cap N_2}$. Hence further assume that
for each $i\in\{1,2\}$, $u$ and $v$ are conjugated by $h_i\in F$
in $G_i$. Therefore the word $h_i^{-1}uh_iv^{-1}$ is equal to the
identity in $G_i$.

By Condition $1.1$ of Theorem $1$ the word $h_1^{-1}uh_1v^{-1}$
can be effectively represented with defining relations $R_1$ of $G_1$ in the
form $\prod_{s=1}^{m_{1}}g_{1,s}r_{1,s}g_{1,s}^{-1}$, where
$r_{1,s}\in R_1$, $g_{1,s}\in F$. By this representation construct
a planar picture $P_{1}$ over the presentation $G_1=\langle A\,
\mid R_1\rangle$ with the boundary label equal to
$h_1^{-1}uh_1v^{-1}$ so that  the edges of $P_{1}$ are labelled by
letters. In addition on the boundary $\partial P_{1}$ of $P_{1}$
fix four points $a_{1}, b_{1}, c_{1}, d_{1}$ not belonging to any
edge and dividing $\partial P_{1}$ into four subpaths so that the
labels of the subpaths $[a_{1}, b_{1}], [b_{1},c_{1}],
[c_{1},d_{1}], [d_{1},a_{1}]$ are identically equal to $h_1^{-1}$,
$u$, $h_1$, $v^{-1}$ respectively. Pasting together the subpaths
$[a_{1}, b_{1}]$ and $[d_{1},c_{1}]$ of $P_{1}$, we obtain an
annulus picture $\overline{P}_{1}$ with the two boundary circles
formed by $[b_{1},c_{1}]$ with the label $u$ and $[d_{1},a_{1}]$
with the label $v^{-1}$. The pasted points $b_{1},c_{1}$ (resp.,
$d_{1},a_{1}$) give a point $({bc})_1$ (resp., $({ad})_1$). The
pasted subpaths $[a_{1}, b_{1}]$ and $[d_{1},c_{1}]$ form a
subpath $Conj_1$.

Similarly changing the index $1$ by the index $2$ in the notation,
by the word $h_2^{-1}u^{-1}h_2v$ construct an annulus picture
$\overline{P}_{2}$ over the presentation $G_2=\langle A\, \mid
R_2\rangle$ with the two boundary circles formed by $[b_{2},
c_{2}]$ with the label $u^{-1}$ and $[d_{2}, a_{2}]$ with the
label $v$.

Pasting together $\overline{P}_{1}$ over $\overline{P}_{2}$ by
their boundaries we obtain a picture $P$ on the torus $T$ over the
presentation $G=\langle A\, \mid R_1\cup R_2 \rangle$. The pasted
circles $[d_{1},a_{1}]$ and $[a_{2}, d_{2}]$ (resp.,
$[b_{1},c_{1}]$ and $[c_{2}, b_{2}]$)  form a circle
$\underline{Equ}$ (resp., $\overline{Equ}$). The pasted points
$({ad})_1$ and $({ad})_2$ ($({bc})_1$ and $({bc})_2$) form a point
$p_v$ ($p_u$). By $Conj$ denote a circle formed by the pasted
subpaths $Conj_1$ and $Conj_2$. The circles $\underline{Equ}$ and
$\overline{Equ}$ will be called {\it the equators.} The points
$p_u$, $p_v$ will be called {\it the poles}.

So $Lab_{p_v}(\underline{Equ})$ is equal to $v^{-1}$ or $v$,
$Lab_{p_u}(\overline{Equ})$ is equal to $u$ or $u^{-1}$ depending
on the direction of travelling along $\underline{Equ}$ and
$\overline{Equ}$. Fix the positive direction of travelling along
the equators so that $Lab^{+}_{p_v}(\underline{Equ})$ is equal to
$v$, $Lab^{+}_{p_u}(\overline{Equ})$ is equal to $u$.

The equators $\underline{Equ}$ and $\overline{Equ}$ divide the
torus $T$ into two annulus (corresponding to $\overline{P}_{1}$
and $\overline{P}_{2}$). The annulus containing the vertices with
labels from $R_1$ (resp., $R_2$) will be called {\it the
$R_1$-annulus} (resp., {\it the $R_2$-annulus}).

In the sequel we will use {\it admissible moves} to transform the
picture $P$ on $T$.  A move is called {\it admissible} if after
the move,

(i) $Lab^{+}_{p_v}(\underline{Equ})$ (resp.,
$Lab^{+}_{p_u}(\overline{Equ})$) is replaced by a word equal to
$v$ (resp., $u$) to within elements from $N_1\cap N_2$;

(ii) all generalized vertices
with labels from $N_1$ are only in the $R_1$-annulus, all generalized
vertices with labels from $N_2$ are only in the $R_2$-annulus, where a {\it generalized vertex} is a vertex (one can consider it as a 'small' planar picture)
with the label equal to an arbitrary (not-necessary reduced) word of $N_1$ (or $N_2$);

(iii) the equators and $Conj$ remain unchanged.

\begin{pro}
Let the presentation $G=\langle A\, \mid R_1\cup R_2 \rangle$ be
atorical (Condition 2.2 of Theorem $1$). Then there exists a
finite sequence of admissible moves
of $P$, at the end of which the labels of the equators will have
one of the following form:

(1) $Lab^{+}_{p_u}(\overline{Equ}) =
\alpha(\omega'\nu_{1}\nu_{2})\alpha^{-1}$,
$Lab^{+}_{p_v}(\underline{Equ}) =
\beta(\omega'\nu_{1}\nu_{2})\beta^{-1}$;

(2) $Lab^{+}_{p_u}(\overline{Equ}) =
\alpha(\omega^{k}\nu_{1}\omega^{-l}\nu_{2}\omega^l)\alpha^{-1}$,
$Lab^{+}_{p_v}(\underline{Equ}) =
\beta(\omega^{k}\nu_{1}\nu_{2})\beta^{-1}$;

\noindent where $\nu_{i}\in N_i$, $\alpha,\beta,\omega,\omega'\in
F$, $l, k\in \mathbb{Z}$, $l\neq 0$ can be determined by $P$ at
the end.
\end{pro}

To prove Theorem $1$ let us use Assertion $1$, which will be
proved in Section $\ref{predl1}$. By Assertion $1$ we have two
possibilities for representation of $u$ and $v$ to within elements
from $N_1\cap N_2$. If $u$ and $v$ have the form (1), they are
evidently conjugated in $F/{N_1\cap N_2}$ by the word $h =
\alpha^{-1}\beta$. Consider the case, when $u$ and $v$ have the
form (2). The following notations will be used: $$\widetilde{u} =
\alpha^{-1}u \alpha, \,\,\,\widetilde{v} = \beta^{-1}v\beta;$$
 $$Roots_{G_1}(\widetilde{v}) = \{c\in F\,\,\, |\,\,\, \exists s = s(c)\in
\mathbb{Z} : \widetilde{v} = c^s\,\,\, {\text в }\,\,\, G_1\};$$
 $$Roots_{G_2}(\widetilde{v}) = \{d\in F\,\,\, |\,\,\, \exists t = t(d)\in
\mathbb{Z} : \widetilde{v} = d^t\,\,\, {\text в }\,\,\, G_2\}.$$

\begin{lem}
Let the presentation $G=\langle A\, \mid R_1\cup R_2 \rangle$ be
atorical (Condition 2.2 of Theorem $1$) and $u =
\alpha(\omega^{k}\nu_{1}\omega^{-l}\nu_{2}\omega^l)\alpha^{-1}$,
$v = \beta(\omega^{k}\nu_{1}\nu_{2})\beta^{-1}$ in $F/{N_1\cap
N_2}$. Then $u$ and $v$ are conjugated in $F/{N_1\cap N_2}$ if and
only if there exist  $c\in Roots_{G_1}(\widetilde{v})$, $d\in Roots_{G_2}(\widetilde{v})$, $\bar{s},
\bar{t}\in \mathbb{Z}$ with $0\leq\bar{s} < s(c), 0\leq\bar{t} < t(d)$ such
that $d^{-\bar{t}}c^{-\bar{s}}\omega^l$ belongs to the
cyclic subgroup $\langle \widetilde{v} \rangle$ of $G$.
\end{lem}
\noindent{\it Proof of Lemma 2.} Assume that there exists a word
$h\in F$ such that the equality $u = h^{-1}vh$ holds in
$F/{N_1\cap N_2}$. Then the equality $u = h^{-1}vh$ holds both in
$F/{N_1}$ and $F/{N_2}$. It is clear that $u$ and $v$ are
conjugated by $h$ if and only if $\widetilde{u}$ and
$\widetilde{v}$ are conjugated by the word $x$, where $x =
\alpha^{-1}h\beta$. Hence further we will consider $\widetilde{u}$
and $\widetilde{v}$ and investigate $x$.

In $G_1 = F/{N_1}$, $\widetilde{u} =
\omega^{-l}(\omega^{k}\nu_{2})\omega^l$ and $\widetilde{v} =
\omega^{k}\nu_{2}$. Since $\widetilde{u} = x^{-1}\widetilde{v}x$
in $G_1$, we have that ${\omega^l}x^{-1}$ and $\widetilde{v}$
commutate in $G_1$. By Condition 2.2 of Theorem $1$ the
presentation $G = \langle A\, \mid R_1\cup R_2 \rangle$ is
atorical, hence, the presentation $G_1 = \langle A\, \mid R_1
\rangle$ is also atorical. By Theorem 13.5 \cite{olsh} it follows
that  there exists $c\in Roots_{G_1}(\widetilde{v})$ such that
$\widetilde{v} = c^s$, ${\omega^l}x^{-1} = c^{m_1}$ in $G_1$ for
some $s=s(c), m_1\in \mathbb{Z}$. On the other hand, $\widetilde{u}$ and
$\widetilde{v}$ are equal to $\omega^{k}\nu_{1}$ in $G_2 =
F/{N_2}$. Since $\widetilde{u} = x^{-1}\widetilde{v}x$ in $G_2$,
we have that $x^{-1}$ and $\widetilde{v}$ commutate in $G_2$.
By Condition 2.2 of Theorem $1$ and Theorem 13.5 \cite{olsh}
there exists $d\in Roots_{G_2}(\widetilde{v})$ such that $\widetilde{v} =
d^t$, $x^{-1} = d^{m_2}$ in $G_2$ for some $t=t(d), m_2\in \mathbb{Z}$.

It follows from the equalities ${\omega^l}x^{-1} = c^{m_1}$ in
$G_1$ and $x^{-1} = d^{m_2}$ in $G_2$ that $\omega^l =
c^{m_1}d^{-m_2}$ in $G = F/{N_1N_2}$. Since $\widetilde{v} = c^s$
in $G_1$, $\widetilde{v} = d^t$ in $G_2$, we have $\omega^l =
c^{\bar{s}}d^{\bar{t}}\widetilde{v}^p$ in $G$ for $0\leq \bar{s} <
s, 0\leq \bar{t} < t$ and some integer $p$, that is,
$d^{-\bar{t}}c^{-\bar{s}}\omega^l = \widetilde{v}^p$ in $G$.

Conversely, suppose $d^{-\bar{t}}c^{-\bar{s}}\omega^l =
\widetilde{v}^p$ in $G = F/{N_1N_2}$. Let us prove that $u$ and
$v$ are conjugated in $F/{N_1\cap N_2}$. Since
$d^{-\bar{t}}c^{-\bar{s}}\omega^l = \widetilde{v}^p$ in $G$, the
word $\widetilde{v}^{-p}d^{-\bar{t}}c^{-\bar{s}}\omega^l$ is
represented in the form $\widetilde{\nu_2}\widetilde{\nu_1}^{-1}$
for some words $\widetilde{\nu_{i}}\in N_i$ ($i=1,2$). Therefore
we have the equality $c^{-\bar{s}}\omega^l\widetilde{\nu_1} =
d^{\bar{t}}\widetilde{v}^{p}\widetilde{\nu_2}$ in $F$. Let us
verify that we can take $c^{-\bar{s}}\omega^l\widetilde{\nu_1} =
d^{\bar{t}}\widetilde{v}^{p}\widetilde{\nu_2}$ as $x$. Indeed, in
$G_1$ we have $$x^{-1}\widetilde{v}x = x^{-1}c^sx =
\omega^{-l}c^{\bar{s}}c^sc^{-\bar{s}}\omega^l =
\omega^{-l}c^s\omega^l = \omega^{-l}\widetilde{v}\omega^l =
\omega^{-l}(\omega^{k}\nu_{2})\omega^l = \widetilde{u}.$$ In $G_2$
we have
$$x^{-1}\widetilde{v}x = x^{-1}d^tx =
\widetilde{v}^{-p}d^{-\bar{t}}d^td^{\bar{t}}\widetilde{v}^{p} =
\widetilde{v}^{-p}\widetilde{v}\widetilde{v}^{p} = \widetilde{v} =
\widetilde{u}.$$  Hence, $x^{-1}\widetilde{v}x = \widetilde{u}$ in
$F/{N_1\cap N_2}$. Therefore $u = h^{-1}vh$ in $F/{N_1\cap N_2}$
for $h = {\alpha}x\beta^{-1}$. \edvo

By Lemma 2 we get the following algorithm.

By the word $\widetilde{v}$ determine finite sets $Roots_{G_1}(\widetilde{v})$, $Roots_{G_2}(\widetilde{v})$ (it is possible by Condition 1.2 of Theorem $1$).  For each $c\in Roots_{G_1}(\widetilde{v})$ and $d\in Roots_{G_2}(\widetilde{v})$, using Condition 1.1 of Theorem $1$, find the numbers $s=s(c), t=t(d)\in \mathbb{Z}$ with the least absolute values such that $\widetilde{v} = c^s$ in $G_1$
and $\widetilde{v} = d^t$ in $G_2$.
Using Condition 2.1 of Theorem $1$, verify whether there exists an integer $p$ such that
$d^{-\bar{t}}c^{-\bar{s}}\omega^l = \widetilde{v}^p$ in $G =
F/{N_1N_2}$ for some integers $\bar{s}, \bar{t}$ with
$0\leq\bar{s} < s(c), 0\leq\bar{t} < t(d)$.
If such $p$ is found,
express $\widetilde{v}^{-p}d^{-\bar{t}}c^{-\bar{s}}\omega^l$ with
defining relations $R_1\cup R_2$ of $G$ (it is possible by
Condition 2.1 of Theorem $1$) and represent
$\widetilde{v}^{-p}d^{-\bar{t}}c^{-\bar{s}}\omega^l$ in the form
$\widetilde{\nu_2}\widetilde{\nu_1}^{-1}$, where
$\widetilde{\nu_{i}}\in N_i$ ($i=1,2$). One can take
${\alpha}c^{-\bar{s}}\omega^l\widetilde{\nu_1}\beta^{-1}$ as a
word $h$ conjugating $u$ and $v$. If for any $c\in Roots_{G_1}(\widetilde{v})$ and $d\in Roots_{G_2}(\widetilde{v})$
there is no such $p$, conclude that $u$ and $v$ are not conjugated
in $F/{N_1\cap N_2}$. So Theorem $1$ is proved. \edvo

\section{Admissible moves using in the proof of Assertion 1.}
Below any domain $M\subset T$ homeomorphic to the square $\{
(x,y)\in \mathbb R^2 \mid \, -1<x<1, -1<y<1\,  \} $ together with
vertices and parts of edges belonging to $M$ will be called {\it a
map}. For a given path (an edge) on the torus $T$,  any part of
the path (the edge) homeomorphic to $\{ x\in \mathbb R \mid \, -1\leq x\leq
1\, \}$ will be called {\it a segment} of the path (of the edge).
We will say that a domain on the torus {\it contains nothing},
if it does not contain poles, vertices and segments of edges of $P$. We will say that a domain on the torus {\it contains absolutely nothing}, if it contains nothing and there is no point from $\underline{Equ}\cup \overline{Equ}\cup Conj$ in it.

\noindent{\it 1) Isotopy.}
\newline
{\it An isotopy} of the picture $P$ is defined by replacing $P$ by
a picture $F_1(P)$, where $F_t :T\times [0,1] \to T\times [0,1]$
is a continuous isotopy of the torus $T$ such that
\begin{itemize}
\item[(i)] $F_t$ leaves fixed all vertices and the both poles, i.e.
for each $t\in [0,1]$ and each vertex $V_i$, $F_t(V_i)=V_i$,
$F_t(p_u)=p_u$, $F_t(p_v)=p_v$;
\item[(ii)] for each $t\in [0,1]$ and each edge $E_j$ the
intersection of $F_t(E_j)$ and $\underline{Equ}$,
$\overline{Equ}$, $Conj$ consists of a finite number of points,
moreover, if $\underline{Equ}$, or
$\overline{Equ}$, or $Conj$ intersects $F_1(E_j)$, then it crosses it, and doesn't just touch it.
\end{itemize}

\unitlength=1mm
\begin{picture}(120,50)(13,-5)
       \thicklines
       \multiput(30,10)(55,0){2}{\line(1,0){30}}
       \put(70,10){\vector(1,0){5}}
       \thinlines
       \multiput(21,10)(96,0){2}{$\underline{Equ}$}
       \multiput(40,25)(55,15){2}{\vector(0,-1){15}}
       \multiput(50,10)(55,15){2}{\vector(0,1){15}}
       \multiput(45,10)(55,15){2}{\oval(10,20)[b]}
       \multiput(35,12)(55,15){2}{\rm g}
       \multiput(53,22)(55,15){2}{\rm g}

\end{picture}
\begin{center} Fig. 1
\end{center}

An isotopy of $P$ is an admissible move because either it
corresponds to a succession of free insertions or free deletions
in $Lab^{+}_{p_v}(\underline{Equ})$ and
$Lab^{+}_{p_u}(\overline{Equ})$ or it does not change
$Lab^{+}_{p_v}(\underline{Equ})$ and
$Lab^{+}_{p_u}(\overline{Equ})$ at all (see Fig.1).

\noindent{\it $2)$ Deletion of a superfluous loop (this is a
particular case of isotopy).}
\newline
Let $\underline{Equ}$ (resp., $\overline{Equ}$, $Conj$) intersect
any edge $E$ in two points, which divide $\underline{Equ}$ (resp.,
$\overline{Equ}$, $Conj$) into two parts so that one of these
parts $\zeta$ does not intersect any edge and does not contain the
poles. By $\vartheta$ denote the segment of $E$ between
these points. If a disk on the torus $T$ encircled by the circle
$\zeta\sqcup\vartheta$ contains absolutely nothing inside, then $\vartheta$
is called {\it a superfluous loop}. It is clear that superfluous
loops do not contribute to the corresponding equatorial label
(considered as an element of the free group). Therefore
superfluous loops can be removed (see Fig.1).

\noindent{\it $3)$ Bridge moves.}
\newline
Assume that a map $M$ contains absolutely nothing except for two segments of edges $\{ x=-1/2, -1<y<1
\}$ {and} $\{ x=1/2, -1<y<1 \},$ which are contrariwise oriented
and labelled by the same word $g$. A transformation of $P$ is
called {\it a bridge move} if it does not change $P$ out of $M$
and change $P$ inside $M$ as is shown on Fig 2. A bridge move is
an admissible move because it does not change the equatorial
labels.

\unitlength=1mm
\begin{picture}(120,50)(13,0)
       \thicklines
       \multiput(30,10)(10,0){3}{\line(1,0){5}}
       \multiput(85,10)(10,0){3}{\line(1,0){5}}
       \multiput(35,40)(10,0){3}{\line(1,0){5}}
       \multiput(90,40)(10,0){3}{\line(1,0){5}}

       \multiput(30,15)(0,10){3}{\line(0,1){5}}
       \multiput(85,15)(0,10){3}{\line(0,1){5}}
       \multiput(60,10)(0,10){3}{\line(0,1){5}}
       \multiput(115,10)(0,10){3}{\line(0,1){5}}

       \put(70,25){\vector(1,0){5}}
       \thinlines
       \put(40,40){\vector(0,-1){30}}
       \put(50,10){\vector(0,1){30}}
       \multiput(95,17)(0,23){2}{\vector(0,-1){7}}
       \multiput(105,10)(0,23){2}{\vector(0,1){7}}
       \put(100,33){\oval(10,10)[b]}
       \put(100,17){\oval(10,10)[t]}
       \multiput(35,12)(55,0){2}{\rm g}
       \multiput(53,37)(55,0){2}{\rm g}
       \multiput(90,35)(18,-21){2}{\rm g}
\end{picture}

\begin{center}
Fig. 2
\end{center}

\noindent$4)$ {\it Uniting of edges.}
\newline
Let $E_1$ and $E_2$ be two edges-circles with labels $g_1$ and
$g_2$, which side by side intersect $\underline{Equ}$,
$\overline{Equ}$ and $Conj$ and bound on the torus an
annulus, containing nothing, or $E_1$ and $E_2$ be two edges with labels $g_1$ and
$g_2$, which join the same vertices, side by side intersect
$\underline{Equ}$, $\overline{Equ}$ and $Conj$ and encircle on the torus a disk, containing nothing. Remove $E_2$.  If $E_1$ and $E_2$ had the
same orientation, label $E_1$ by $g_1g_2$ or $g_2g_1$, otherwise
label $E_1$ by $g_1g_2^{-1}$ or $g_2^{-1}g_1$. The label for $E_1$
should be chosen so that the contribution of this label to the
equatorial labels remains the same as the contribution of the both
edges $E_1$ and $E_2$. We will assume that the multiplication of $g_1$ and $g_2^{\pm 1}$ is free.

\noindent{\it $5)$ Cutting of complete dipoles.}
\newline
{\it A complete $R_1$-dipole} is a dipole $D_1$ such that the
labels  of its vertices is equal to $r_1^{\pm 1}\in R_1\setminus
R_2$ and its vertices are joined by a single edge $E_1$ with the
label $r_1$.

Consider a map $M$ in the $R_1$-annulus such that $M$
contains absolutely nothing except for a segment of $E_1$: $\{ x=0 \}$, starting
at the point $(0,-1)$ and ending at the point $(0,1)$. Cut out $M$
from $P$ and paste a new map $M^{\prime}$ instead of $M$. The new
map $M^{\prime}$ contains absolutely nothing except for two vertices $V'$,
$V''$ with the labels $r_1$, $r_1^{-1}$ and two edges one of which
starts at $(0,-1)$ and ends at $V'$, and the other one starts at
$V''$ and ends at $(0,1)$. As a result one has two complete
$R_1$-dipoles instead of one. This move is admissible because it
does not change the equatorial labels.

Similarly one can define {\it a complete $R_2$-dipole} whose
vertices are labelled by $r_2^{\pm 1}\in R_2\setminus R_1$ and a
corresponding move performed in the $R_2$-annulus. Similarly one
can define {\it a complete mixed dipole} whose vertices are
labelled by $r^{\pm 1}\in R_1\cap R_2$.

\noindent{\it $6)$ Conjugation of dipoles.}
\newline
Let $n_1\in N_1$. {\it A generalized $N_1$-dipole} with the label
$n_1$ is two generalized vertices with the labels $n_1^{\pm 1}$
and a single edge with the label $n_1$, joining them. For example,
a complete $R_1$-dipole is generalized one labelled by $r_1\in
R_1\subset N_1$.

Let $D_1$ be a generalized $N_1$-dipole with the label $n_1$ and
$C$ be an edge-circle with a label $f\in F$, encircling on $T$ a
disk, containing nothing except for $D_1$. In addition the edge of
$D_1$ and the edge-circle $C$ side by side intersect
$\underline{Equ}$, $\overline{Equ}$ and $Conj$ and contribute
$(fn_1f^{-1})^{\pm 1}\in N_1$ to the labels of $\underline{Equ}$,
$\overline{Equ}$ and $Conj$. Remove $C$ and label the edge of
$D_1$ by $fn_1f^{-1}$ and its generalized vertices by
$(fn_1f^{-1})^{\pm 1}$. This move does not change the equatorial
labels, hence it is admissible.

Similarly one can define {\it a generalized $N_2$-dipole} and a
corresponding move of it.

\noindent$7)$ {\it Deletion of a dipoles and an edge-circles not
intersecting the equators.}
\newline
If a generalized dipole or an edge-circle does not intersect
$\underline{Equ}$ and $\overline{Equ}$, then it does not
contribute to $Lab^{+}_{p_v}(\underline{Equ})$ and
$Lab^{+}_{p_u}(\overline{Equ})$. Hence remove it.

\noindent$8)$ {\it Conjugation of a pole.}
\newline
Consider the pole $p_v$ (everything is similar for $p_u$). Let $C$
be an edge-circle with a label $g\in F$ and $C$ encircle on $T$ a
disk containing absolutely nothing except for $p_v$, only one segment of
$\underline{Equ}$ and only one segment of $Conj$. The union of $C$
and $p_v$ is called {\it a conjugated pole $p_v$}. The pole $p_v$
itself will be considered as a conjugated pole (encircled by an
edge-circle $C$ with the label equal to the identity of the free
group). If a conjugated pole $p_v$ is surrounded
in the same way by an edge-circle $\tilde{C}$, then unite $C$ and $\tilde{C}$.
This move does not change the equatorial labels, hence it is
admissible.

\noindent{\it $9)$ Deletion of a one-sided dipole.}
\newline
Let the edge of a generalized $N_1$-dipole $D_1$ (everything is
similar for a generalized $N_2$-dipole) with the label $n_1\in
N_1$ do not intersect $Conj$ and intersect only one of the
equators (for definiteness, $\overline{Equ}$) and only at two
points. Then $D_1$ is called {\it a one-sided $N_1$-dipole}.

There exists a closed disk $O$ containing absolutely nothing except for $D_1$
and two segments $[s_1, s_2]$ and $[t_1,t_2]$ of $\overline{Equ}$,
where the points $s_1, s_2, t_1,t_2$ belong to ${\partial
O\cap\overline{Equ}}$. Note that the labels of $[s_1, s_2]$ and
$[t_1,t_2]$ are equal to $n_1$ and $n_1^{-1}$ respectively, i.e.,
to the labels of $D_1$. In addition either $[s_2, t_1]$ or $[t_2,
s_1]$ does not contain the pole. For definiteness let us assume
that it is $[s_2, t_1]$. The points $s_2$, $t_1$ divide $\partial
O$ into two segments. By $\varrho$ denote such of them which
contains no points of the $R_1$-annulus. Then the closed path
$[s_2, t_1]\cup\varrho$ encircles a planar picture over the
presentation $G=\langle A\, \mid R_2 \rangle$. By Lemma $1$ the
label $n_2$ of $[s_2, t_1]\cup\varrho$ belongs to $N_2$. Since no
edges intersect $\varrho$, $n_2$ is the label of $[s_2, t_1]$. So
the label of $[s_1, s_2]\cup [s_2, t_1]\cup [t_1,t_2]$ is equal to
$n_1n_2n_1^{-1}$. Remove $D_1$ from $P$. The label of $[s_1,
s_2]\cup [s_2, t_1]\cup [t_1,t_2]$ becomes equal to $n_2$. This
move does not change $Lab^{+}_{p_u}(\overline{Equ})$ to within
$n_1n_2n_1^{-1}{n_2}^{-1}\in N_1\cap N_2$. Hence this move is
admissible.

\noindent{\it $10)$ Permutation of two-sided dipoles.}
\newline
Let the edge of a generalized $N_1$-dipole $D_1$ with the label
$n_1\in N_1$ do not intersect $Conj$ and intersect each of the
equators $\overline{Equ}$ and $\underline{Equ}$ exactly at one
point. Then $D_1$ is called {\it a two-sided $N_1$-dipole}.

There is an open disk $O_1$ containing absolutely nothing except for $D_1$
and two segments $[s_1, t_1]\in\overline{Equ}$ and
$[q_1,p_1]\in\underline{Equ}$, where the points $s_1, t_1$ belong
to ${\partial O_1\cap\overline{Equ}}$, the points $q_1,p_1$ belong
to ${\partial O_1\cap\underline{Equ}}$. Note that the labels of
$[s_1, t_1]$ and $[q_1,p_1]$ are equal to $n_1$ and $n_1^{-1}$
respectively, i.e., to the labels of $D_1$.

Similarly one can define {\it a two-sided $N_2$-dipole}.
Substituting $2$ instead of $1$ in the above notations for the
two-sided $N_1$-dipole, one gets the same notations for a two-sided
$N_2$-dipole.

Now let both a two-sided  $N_1$-dipole $D_1$ and a two-sided
$N_2$-dipole $D_2$ be in $P$. The points $s_1, s_2, t_1, t_2$
divide $\overline{Equ}$ into four segments. Assume that one of
them (say $[t_1, s_2]$) does not intersect any edge and does not
contain the pole. Then the label of the segment $\sigma=[s_1,
t_1]\cup [t_1, s_2]\cup [s_2, t_2]$ is equal to $n_1n_2$. Permute
the segments $[s_1, t_1]$ and $[s_2, t_2]$ (see Fig. 3).

\unitlength=1mm
\begin{picture}(120,50)(13,0)
       \thicklines
       \multiput(30,10)(10,0){3}{\line(1,0){5}}
       \multiput(85,10)(10,0){3}{\line(1,0){5}}
       \multiput(35,40)(10,0){3}{\line(1,0){5}}
       \multiput(90,40)(10,0){3}{\line(1,0){5}}

       \multiput(30,15)(0,10){3}{\line(0,1){5}}
       \multiput(85,15)(0,10){3}{\line(0,1){5}}
       \multiput(60,10)(0,10){3}{\line(0,1){5}}
       \multiput(115,10)(0,10){3}{\line(0,1){5}}

       \put(70,25){\vector(1,0){5}}

        \thicklines
        \multiput(30,25)(55,0){2}{\line(1,0){30}}
        \thinlines
        \multiput(20,25)(97,0){2}{$\overline{Equ}$}

        \put(38,10){\vector(0,1){20}}
        \put(38,31){\circle{2}}
        \qbezier[10](35,30)(35,34)(38,34)
        \qbezier[10](41,30)(41,34)(38,34)
        \qbezier[30](35,10)(35,20)(35,30)
        \qbezier[30](41,10)(41,20)(41,30)

        \put(52,40){\vector(0,-1){20}}
        \put(52,19){\circle*{2}}
        \qbezier[10](49,20)(49,16)(52,16)
        \qbezier[10](55,20)(55,16)(52,16)
        \qbezier[30](49,40)(49,30)(49,20)
        \qbezier[30](55,40)(55,30)(55,20)
        \put(32,26){$\rm s_1$}
        \put(42,26){$\rm t_1$}
        \put(46,22){$\rm s_2$}
        \put(56,22){$\rm t_2$}

        \put(93,25){\vector(0,-1){5}}
        \put(93,19){\circle*{2}}
        \qbezier[10](90,20)(90,16)(93,16)
        \qbezier[10](96,20)(96,16)(93,16)
        \qbezier[10](90,20)(90,22)(90,25)
        \qbezier[10](96,20)(96,22)(96,25)

        \put(93,10){\vector(0,1){2}}
        \qbezier[30](94,15)(105,17)(104,25)
        \qbezier[10](90,10)(89,14)(94,15)
        \qbezier(96,14)(107,14)(107,25)
        \qbezier(93,10)(92,13)(96,14)
        \qbezier[30](99,12)(110,13)(110,25)
        \qbezier[5](96,10)(95,12)(99,12)

        \put(107,25){\vector(0,1){5}}
        \put(107,31){\circle{2}}
        \qbezier[10](104,30)(104,34)(107,34)
        \qbezier[10](110,30)(110,34)(107,34)
        \qbezier[10](104,25)(104,28)(104,30)
        \qbezier[10](110,25)(110,28)(110,30)

        \put(107,40){\vector(0,-1){2}}

        \qbezier[30](106,35)(95,33)(96,25)
        \qbezier[10](110,40)(111,36)(106,35)
        \qbezier(104,36)(93,36)(93,25)
        \qbezier(107,40)(108,37)(104,36)
        \qbezier[30](101,38)(90,37)(90,25)
        \qbezier[5](104,40)(105,38)(101,38)

        \put(87,22){$\rm s_2$}
        \put(97,22){$\rm t_2$}
        \put(101,26){$\rm s_1$}
        \put(111,26){$\rm t_1$}
\end{picture}

\begin{center}
Fig. 3
\end{center}

After this move, the label of the new segment $\sigma$ become equal
to $n_2n_1$. This move is admissible, since after it
$Lab^{+}_{p_u}(\overline{Equ})$ is not changed to within  the word
$n_2^{-1}n_1^{-1}n_2n_1\in [N_1, N_2].$

Similarly one can define the same move for $\underline{Equ}$.

\noindent{\it $11)$ Moving of an edge over a dipole or a pole.}
\newline
Let $X$ be a generalized $N_1$- or $N_2$-dipole (resp., a
conjugated $p_u$ or $p_v$ pole), $O_1, O_2, O_3$ be three closed
disks on the torus $T$ containing nothing except for $X$ such that $O_3\subset
O_2\backslash {\partial O_2}, O_2\subset O_1\backslash {\partial
O_1}$. Let $E$ be an edge with a label $g\in F$ such that $E\cap
O_1 = \emptyset$ and there exists a simple path $\gamma$ joining
points $o\in E$ and $o_3\in
\partial O_3$, intersecting $E$ and $\partial O_1, \partial O_2$ exactly at one
point, not intersecting other edges, the equators and $Conj$ and
not passing through any vertex. Thus $Lab^+(\gamma) = g$. Put in $P$ two
contrariwise oriented edge-circles $C_1 =
\partial O_1$ and $C_2 = \partial O_2$ labelled by $g\in F$  so that
$Lab^+(\gamma)$ becomes identically equal to $gg^{-1}g$. Apply the
bridge move to $E$ and $C_1$, the conjugation to $X$ and $C_2$. It
is clear that this move is admissible, because either it
corresponds to an insertion of  inverse words in
$Lab^{+}_{p_v}(\underline{Equ})$ and
$Lab^{+}_{p_u}(\overline{Equ})$, or it does not change
$Lab^{+}_{p_v}(\underline{Equ})$ and
$Lab^{+}_{p_u}(\overline{Equ})$ at all.

\section{Proof of Assertion 1.}\label{predl1}

{\bf STEP 1.} {\it Extraction of complete dipoles.}

Since the presentation $G=\langle A\, \mid R_1\cup R_2 \rangle$ is
atorical, there exists a dipole $D$ in the picture $P$ on the
torus $T$, i.e., there exists a couple of vertices $V_1$ and $V_2$
with mutually inverse labels $r$ and $r^{-1}$ such that $V_1$ and
$V_2$ are connected by an edge $\rho$. Applying the bridge moves
no more than $|r|-1$ times, we obtain that all edges go from $V_1$
to $V_2$ side by side in a parallel way to $\rho$ and $\rho$
remains unchanged. Unite these edges. This makes the dipole $D$
complete.

Now the picture $P$ consists of two disjoint subpictures
$P_1\sqcup P_2$, one of which (say $P_1$) contains nothing except for the
complete dipole $D$.  The subpicture
$P_2$ is a picture on the torus over presentation $G=\langle A\,
\mid R_1\cup R_2 \rangle$. Besides $P_2$ contains two fewer
vertices than $P$. Repeating the above
procedure for $P_2$, and so on, we will eventually reduce $P$ to
$m_V/2$ complete dipoles and edge-circles, where $m_V$ is the number of vertices in $P$.

{\bf STEP 2.} {\it A move after which dipoles do not intersect
$Conj$.}

After Step 1 the picture $P$ consists of edges-circles, complete
$R_1$-, $R_2$-dipoles and complete mixed dipoles. If the edges of
some complete dipoles do not intersect the equators, remove these
complete dipoles. Also remove complete mixed dipoles from $P$.
This changes the equatorial labels by elements from $R_1\cap R_2 \subset N_1\cap N_2$. Now $P$
contains only complete $R_1$-, $R_2$-dipoles and edges-circles.

\noindent{\bf Operation 1.} Consider a complete $R_1$-dipole $D$
(the case of a complete $R_2$-dipole is similar). Let its edge
intersect $Conj_1$ at points $o_1,...,o_{\tilde{m}}$. Near by $o_i$
($i=1,...,\tilde{m}$) cut $D$ into three complete dipoles $D_1, D_2, D_3$,
one of which ($D_2$) lies in the $R_1$-annulus as a whole and its
edge intersects $Conj_1$ exactly at one point (at $o_i$). Remove
$D_2$ from $P$. Repeating the same procedure to each of $\tilde{m}$
intersections, instead of one complete $R_1$-dipole $D$, we obtain
$\tilde{m}+1$ complete $R_1$-dipoles, neither of which intersects
$Conj_1$.

Apply Operation 1 to each of complete $R_1$- and $R_2$-dipoles.
This gives that the edges of the complete $R_1$-dipoles do not
intersect $Conj_1$ and the edges of the complete $R_2$-dipoles do
not intersect $Conj_2$.

\noindent{\bf Operation 2.} Consider $Conj_1$ (the case of
$Conj_2$ is similar). It can be intersected only by the edges of
complete $R_2$-dipoles and by edges-circles. Let
$\rho_1,...,\rho_{\hat{m}}$ be edges-circles not conjugating the poles and
edges of complete $R_2$-dipoles such that $\rho_1,...,\rho_{\hat{m}}$
intersect $Conj_1$, and we encounter them in the order
$\rho_1,...,\rho_{\hat{m}}$ if we start at the conjugated pole $p_u$ and
travel along $Conj_1$ to the conjugated pole $p_v$. Starting with
$\rho_1$, move consecutively each edge $\rho_i$ over the
conjugated pole $p_u$. This gives that the edges of the complete
$R_2$-dipoles intersect only $Conj_2$. Apply Operation 1 to these
complete $R_2$-dipoles.

After Operation 2 applying to $Conj_1$ and $Conj_2$,  the picture
$P$ consists of edges-circles and only of complete $R_1$- and
$R_2$-dipoles $D_1$, $....,$ $D_{\bar{m}}$ not intersecting $Conj$. For
each $D_i$, by $m_i$ denote the number of intersections of the
equators and the edge of $D_i$. Note that $m_i$ is even. One can
assume that $m_i > 0$, otherwise remove $D_i$ from $P$. If
$m_i>2$, cut $D_i$ into $m_i/2$ complete dipoles each of which
intersects the equators exactly at two points. This move applying
to each $D_i$ makes all dipoles either one-sided or two-sided.
Remove all one-sided dipoles from $P$.

{\bf STEP 3.} {\it Getting rid of contractible edges-circles.}

After Step 2 the picture $P$ consists just of two-sided dipoles
and edges-circles. Call an edge-circle {\it contractible}, if it
divides the torus into two parts one of which is homeomorphic to a
disk. This part will be called the {\it interior} of the
edge-circle.

After the deletions of superfluous loops from the equators and
$Conj$ each contractible edge-circle $C$ belongs to one of the
following types.

I) The interior of $C$ contains absolutely nothing.

II) The interior of $C$ contains nothing except for just one two-sided dipole or just one conjugated pole.

III) In the interior of $C$, there are at least two two-sided dipoles, or at least one
two-sided dipole and at least one conjugated pole, or the both
conjugated poles.

Remove all edges-circles of Type I from $P$. Apply the conjugation
of dipoles or the conjugation of poles to all edges-circles of
Type II. Now just $\check{m}$ edges-circles of Type III remain in $P$.

Call an edge-circle of Type III {\it minimal}, if there is no
other edges-circles of Type III in its interior.

\noindent{\bf Operation 3.} Let $C$ be a minimal edge-circle of
Type III with $\check{m}_1$ two-sided dipoles and $\check{m}_2$ conjugated poles
in its interior, $\check{m}_1 + \check{m}_2 \geq 2$. It is clear that by the
isotopy and the $\check{m}_1 + \check{m}_2 - 1$ bridge moves, $C$ can be reduced
to $\check{m}_1 + \check{m}_2$ edges-circles of Type II. Apply the conjugation to
each of these edges-circles of Type II.

Operation 3 gives a picture $P$ with one fewer edges-circles.
Hence after no more than $\check{m}$ applications of Operation 3, the
picture $P$ will contain just non-contractible edge-circles and
two-sided dipoles.

{\bf STEP 4.} {\it Uniting of non-contractible edges-circles.}

After Step 3 the picture $P$ contains just non-contractible
edges-circles $Z_1, ..., Z_m$ and two-sided dipoles. Cutting out
one of the edges-circles ($Z_1$) from the torus converts the torus
to a surface $\Omega$ homeomorphic to an annulus. In $\Omega$ any
closed simple not contractible path ($Z_i$, $i\neq 1$) is
homotopic to the boundary ($Z_1$) and to any other closed simple
not contractible path ($Z_j$, $j\neq 1, i$) disjoint with it. The
edges-circles $Z_2, ..., Z_m$ divide $\Omega$ into $m$ disjoint
parts $\Omega_1$, ..., $\Omega_m$ each homeomorphic to an annulus.
Assume that the edges-circles $Z_2, ..., Z_m$ are numbered so that
$\Omega_1$ are bounded by $Z_1$ and $Z_2$, $\Omega_2$ are bounded
by $Z_2$ and $Z_3$,..., $\Omega_m$ are bounded by $Z_m$ and $Z_1$.

Consider $\Omega_1$. If there are conjugated poles or two-sided
dipoles in $\Omega_1$, apply the isotopy and move $Z_2$ over these
dipoles and poles to transpose these poles and dipoles from
$\Omega_1$ to $\Omega_2$, and to approach $Z_2$ and $Z_1$ to each
other so that $Z_1$ and $Z_2$ become parallel and side by side
intersect the equator and $Conj$. Repeat the same procedure for
each $\Omega_i$, $i=2,...,m-1$ to transpose conjugated poles and
generalized dipoles from $\Omega_i$ to $\Omega_{i+1}$. We will
eventually obtain that all conjugated poles and two-sided dipoles
of $P$ are in $\Omega_m$ and $Z_1, ..., Z_m$ are parallel and side
by side intersect the equators and $Conj$. Unite $Z_1, ..., Z_m$.
This gives a single edge-circle $Z$.

{\bf STEP 5.} {\it Disposition of two-sided dipoles in the order.}

Above the orientation on the equators was fixed. If we start at
$p_v$ and travel once around $\underline{Equ}$ in the positive
direction, we encounter a succession of edges of dipoles
$D_1$,...,$D_s$ intersecting $\underline{Equ}$. We say that an
$N_2$-dipole $D_i$ and an $N_1$-dipole $D_j$ form {\it the
inversion} on $\underline{Equ}$, if $i<j$, otherwise they form
{\it the order} on $\underline{Equ}$. In the same way one can
define the inversion and the order on $\overline{Equ}$. We will
say that one circuit along $\underline{Equ}$ (resp.,
$\overline{Equ}$) in the positive direction starting at $p_v$
(resp., $p_u$) is a movement {\it from the left to the right}.

The edge-circle $Z$ is divided by the equators into segments. {\it
Two-sided $R_2$-pieces} (resp., {\it two-sided $R_1$-pieces}) are
such of these segments that do not intersect $Conj$ and lie in the
$R_2$-annulus (resp., in the $R_1$-annulus) at the whole, starting
on one of the equators and ending on the other one.

\begin{lem}
There exists a finite succession of admissible moves that disposes
all edges of two-sided $N_1$-dipoles in the $R_2$-annulus on the
left side of the two-sided $R_2$-pieces of $Z$.
\end{lem}

\noindent{\it Proof of Lemma 3.} If there are no two-sided
$N_1$-dipole or two-sided $R_2$-pieces in $P$, there is nothing to
prove. Otherwise let $m$ be the minimal number of transpositions
to get all edges of two-sided $N_1$-dipoles on the left side of
the two-sided $R_2$-pieces of $Z$. If $m=0$, there is nothing to
prove. Otherwise consider the rightest two-sided $N_1$-dipole
${D}$ which has a two-sided $R_2$-piece $\rho$ on the left such
that there are no other $N_1$-dipoles or two-sided $R_2$-pieces
between ${D}$ and $\rho$. Move $\rho$ over $D$ to the right of
${D}$. This decreases $m$ by $1$. Now use induction on $m$. \edvo

The edges of two-sided $N_1$-dipoles consecutively intersect
$\underline{Equ}$ (resp., $\overline{Equ}$). For a given two-sided
$N_1$-dipole, let $o'$ and $o''$ be two consecutive intersections
of its edge and $\underline{Equ}$ (resp., $\overline{Equ}$). By
Lemma $3$, removing superfluous loops, if necessary, either there
are no intersections with $Z$ between $o'$ and $o''$, or there are
intersections with the edges of two-sided $N_2$-dipoles between
$o'$ and $o''$ and $Z$ intersects $\underline{Equ}$ (resp.,
$\overline{Equ}$) between $o'$ and $o''$, gets into the
$R_2$-annulus, envelops a vertex of at least one of these
two-sided $N_2$-dipoles, turns back to $\underline{Equ}$ (resp.,
$\overline{Equ}$) and returns to the $R_1$-annulus.

\begin{lem}
There exists a finite succession of admissible moves that disposes
all two-sided dipoles of $P$ in the order.
\end{lem}

\noindent{\it Proof of Lemma 4.} By  $m'$ denote the number of inversions on
$\underline{Equ}$, by $m''$ the number of inversions on
$\overline{Equ}$. If $m'+m''=0$, there is nothing to prove. Let $m'+m''>0$.

If $m'>0$, at first consider $\underline{Equ}$. Let $D_1$ and
$D_2$ be two neighboring two-sided $N_1$- and $N_2$-dipoles
forming the inversion on $\underline{Equ}$ such that there are no
other dipoles between them. We can assume that $D_2$ is not
enveloped by $Z$, otherwise move $Z$ over $D_2$. The permutation
of $D_1$ and $D_2$ decreases $m'$ by $1$. Induction on $m'$ gives
that all two-sided dipoles form the order on $\underline{Equ}$.

If $m''>0$, apply the same procedure to $\overline{Equ}$. \edvo

\begin{lem}
There exists a finite succession of admissible moves that disposes
all edges of two-sided $N_2$-dipoles in the $R_2$-annulus on the
left side of two-sided $R_1$-pieces of $Z$.
\end{lem}

The proof of Lemma $5$ is similar to the proof of Lemma $3$.

So all two-sided dipoles of $P$ form the order both on
$\overline{Equ}$ and on $\overline{Equ}$. In addition two-sided
$N_1$-dipoles (resp., $N_2$-dipoles) are near by to each other and
intersect the equators side by side. Replace all these two-sided
$N_1$-dipoles (resp., $N_2$-dipoles) by one two-sided dipole
$\Delta_1$ (resp., $\Delta_2$) with the edge's label $\nu_{1}$
(resp., $\nu_{2}$) equal to the product of the labels of all these
two-sided $N_1$-dipoles (resp., $N_2$-dipoles), i.e., $\nu_{1}$
(resp., $\nu_{2}$) belongs to $N_1$ (resp., $N_2$). If $\nu_{1}$
(resp., $\nu_{2}$) is equal to the identity in $F$, remove the
dipole $\Delta_1$ (resp., $\Delta_2$).

{\bf STEP 6.} {\it Finale.}

The picture $P$ can contain at most one edge-circle $Z$, at
most one two-sided $N_1$-dipole $\Delta_1$ (with the label
$\nu_{1}$), at most one two-sided $N_2$-dipole $\Delta_2$ (with
the label $\nu_{2}$) and two conjugated poles $p_u$ and $p_v$. By
$\alpha$ (resp., $\beta$) denote the label of the edge conjugating
the pole $p_u$ (resp., $p_v$). Below $P$ will be transformed by
isotopy, by moving $Z$ over $\Delta_1$ and $\Delta_2$, by conjugation of
poles. For simplicity of notation the labels of $\Delta_1$ and
$\Delta_2$, the labels of edges conjugating $p_u$ and $p_v$ will
be again denoted by $\nu_{1}$, $\nu_{2}$, $\alpha$, $\beta$.

There are three possibility:

{\bf Case $\bf A$.} {\it There is no $Z$ in $P$.}

We have Case (1) of Assertion $1$, i.e.
$Lab^{+}_{p_u}(\overline{Equ}) =
\alpha(\nu_{1}\nu_{2})\alpha^{-1}$,
$Lab^{+}_{p_v}(\underline{Equ}) =
\beta(\nu_{1}\nu_{2})\beta^{-1}$.

{\bf Case $\bf B$.} {\it There is $Z$ in $P$ and
$Z$ is homotopic to $Conj$.}

By isotopy, moving $Z$ over $\Delta_1$, $\Delta_2$ and the
conjugated poles, dispose $Z$ near by $Conj$ in a parallel way to
$Conj$ so that $Z$ intersects each of the equators exactly at one
point. Thus we have Case (1) of Assertion $1$, i.e.,
$Lab^{+}_{p_u}(\overline{Equ}) =
\alpha(\omega'\nu_{1}\nu_{2})\alpha^{-1}$,
$Lab^{+}_{p_v}(\underline{Equ}) =
\beta(\omega'\nu_{1}\nu_{2})\beta^{-1}$, where $\omega'$ is the
label of $Z$.

{\bf Case $\bf C$.} {\it There is the edge-circle $Z$ in $P$ and
$Z$ is homotopic to a simple closed path circuiting  $Conj$ $|k|$
times and the equators $|l|$ times, where $l,k\in \mathbb{Z}$,
$|l|\geq 1$.}

If there is no dipole $\Delta_2$ in $P$, by isotopy and moving $Z$
over $\Delta_1$ and the conjugated poles, dispose $Z$ in such a way
that the $|k|$ circuits of $Z$ along $Conj$ are near by $Conj$ and
the $|l|$ circuits of $Z$ along the equators are in the
$R_1$-annulus. Thus we have Case (1) of Assertion $1$, i.e.,
$Lab^{+}_{p_u}(\overline{Equ}) =
\alpha(\omega^{k}\nu_{1})\alpha^{-1}$,
$Lab^{+}_{p_v}(\underline{Equ}) =
\beta(\omega^{k}\nu_{1})\beta^{-1}$, where $\omega$ is the label
of $Z$.

It remains to consider the case when there exists $\Delta_2$ in
$P$. By isotopy and moving $Z$ over $\Delta_1$, $\Delta_2$ and the
conjugated poles, dispose $Z$ in such a way that the $|k|$
circuits of $Z$ along $Conj$ are near by $Conj$ and the $|l|$
circuits of $Z$ along the equators start in the $R_1$-annulus,
go in a parallel way to each other to the edge of
$\Delta_2$, envelope its vertex after intersecting
$\overline{Equ}$ and return to the $R_1$-annulus. Thus we have
Case (2) of Assertion $1$: $Lab^{+}_{p_u}(\overline{Equ}) =
\alpha(\omega^{k}\nu_{1}\omega^{-l}\nu_{2}\omega^l)\alpha^{-1}$,
$Lab^{+}_{p_v}(\underline{Equ}) =
\beta(\omega^{k}\nu_{1}\nu_{2})\beta^{-1}$, where $\omega$ is the
label of $Z$.
 \edvo

\begin{rem} It follows from the proof of Assertion $1$ that the integer
$L=L(u,v,R_1,R_2)$ such that $|\alpha|$, $|\beta|$, $|\omega|$,
$|l| \leq L$, can be chosen as $90(|h_1| + |h_2|)(1 + 2l_R +
2l_R^2 + ... + 2l_R^{m_V/2-1})$, where $l_R$ is the length of the
longest word of $R_1\cup R_2$, $m_V$ is the number of vertices in
the initial picture $P$.
\end{rem}


\end{document}